\documentclass[12pt]{article}

\usepackage{subfigure}

\usepackage{latexsym}
\usepackage{graphics}
\usepackage{amsmath}
\usepackage{xspace}
\usepackage{amssymb}
\usepackage{psfrag}
\usepackage{epsfig}
\usepackage{amsthm}
\usepackage{pst-all}

\usepackage{color}

\usepackage{bbm}
\usepackage{mathrsfs}
\usepackage{epsfig}
\usepackage{soul,color}
\usepackage{graphicx}
\usepackage{amsfonts}
\usepackage{amsthm}
\usepackage[figuresright]{rotating}
\usepackage{amssymb}
\usepackage{amsmath}
\usepackage{dcolumn}
\usepackage{physics}
\usepackage{float}
\usepackage{bm}
\usepackage{verbatim}
\usepackage[normalem]{ulem}
\usepackage[colorlinks,linkcolor=blue,anchorcolor=blue,citecolor=blue,urlcolor=blue]{hyperref}
\usepackage{mathtools}

\usepackage{algorithm}
\usepackage{algorithmic}

\definecolor{Red}{rgb}{1,0,0}
\definecolor{Blue}{rgb}{0,0,1}
\definecolor{Olive}{rgb}{0.41,0.55,0.13}
\definecolor{Yarok}{rgb}{0,0.5,0}
\definecolor{Green}{rgb}{0,1,0}
\definecolor{MGreen}{rgb}{0,0.8,0}
\definecolor{DGreen}{rgb}{0,0.55,0}
\definecolor{Yellow}{rgb}{1,1,0}
\definecolor{Cyan}{rgb}{0,1,1}
\definecolor{Magenta}{rgb}{1,0,1}
\definecolor{Orange}{rgb}{1,.5,0}
\definecolor{Violet}{rgb}{.5,0,.5}
\definecolor{Purple}{rgb}{.75,0,.25}
\definecolor{Brown}{rgb}{.75,.5,.25}
\definecolor{Grey}{rgb}{.5,.5,.5}

\newcommand{\pr}{\mathbb{P}}
\newcommand{\E}[1]{\mathbb{E}\!\left[#1\right]}

\newcommand{\G}{\mathbb{G}}

\newcommand{\ignore}[1]{\relax}

\newtheorem{theorem}{Theorem}[section]

\newtheorem{lemma}[theorem]{Lemma}

\newtheorem{prop}[theorem]{Proposition}
\newtheorem{proposition}[theorem]{Proposition}
\newtheorem{coro}[theorem]{Corollary}

\newcommand{\ER}{Erd{\"o}s-R\'{e}nyi }

\newcommand{\dg}[1]{{\small \color{red}[David: #1]}}

\definecolor{Red}{rgb}{1,0,0}
\definecolor{Blue}{rgb}{0,0,1}
\definecolor{Olive}{rgb}{0.41,0.55,0.13}
\definecolor{Green}{rgb}{0,1,0}
\definecolor{MGreen}{rgb}{0,0.8,0}
\definecolor{DGreen}{rgb}{0,0.55,0}
\definecolor{Yellow}{rgb}{1,1,0}
\definecolor{Cyan}{rgb}{0,1,1}
\definecolor{Magenta}{rgb}{1,0,1}
\definecolor{Orange}{rgb}{1,.5,0}
\definecolor{Violet}{rgb}{.5,0,.5}
\definecolor{Purple}{rgb}{.75,0,.25}
\definecolor{Brown}{rgb}{.75,.5,.25}
\definecolor{Grey}{rgb}{.5,.5,.5}
\definecolor{Pink}{rgb}{1,0,1}
\definecolor{DBrown}{rgb}{.5,.34,.16}
\definecolor{Black}{rgb}{0,0,0}

\usepackage{float}

\author{
David Gamarnik \thanks{Operations Research Center,  Center for Statistics and Data Sciences, IDSS, Massachusetts Institute of Technology,
\href{mailto:gamarnik@mit.edu}{gamarnik@mit.edu}}
}

\begin{document}

\title{The satisfiability conjecture for random k-SAT formulas. Human-friendly exposition of the recent OpenAI resolution}
\date{\today}

\maketitle

\begin{abstract}
The satisfiability conjecture was one of the most important open problems in discrete probability. It was resolved recently first 
by Carenini~\cite{GCarenini}, and  shortly after by OpenAI~\cite{k-sat-OpenAI}, among other 370 problems solved by its internal model. 
The conjecture concerned the satisfiability property of  random instances of the k-SAT problem.
It posits the existence of a fixed value for the clause-to-variables ratio around which the probability that the formula is satisfiable drops sharply from nearly one
to nearly zero as the ratio increases. While the case of large constant  $k$  was resolved a decade ago by Ding, Sly and Sun~\cite{ding2015proof,ding2022proof} 
using the heavy machinery
of the second moment method, the case of general $k$ eluded  proof until ~\cite{GCarenini} and~\cite{k-sat-OpenAI}.

In this short note we provide a ''human-friendly'' and detailed exposition of the OpenAI's resolution of the conjecture.   
\end{abstract}



\section{Introduction}

\ignore{
\dg{
\begin{itemize}
\item
Comment on proof not relying on Friedgut. 

\item Comment on Capelini proof

\item Ding-Sly-Sun large k resolution

\item Coloring open problem? Proving validity of the satisfiability threshold?

\item Non-convex problems. Min-bisection
\end{itemize}
}
}

The satisfiability conjecture was one of the most important  open problems in   discrete probability and theoretical computer science areas
in the last several decades.
It concerns random instances of a particular class of constraint satisfaction problems known as k-SAT problem. The k-SAT problem concerns
the satisfiability property of a k-SAT Boolean  formula with $m$ clauses of size $k$ and $n$ variables, see the next section for the formal 
description of the model and the  statement of the conjecture.  
The conjecture posits that for each $k$ there exists a  fixed value $\lambda^*=\lambda^*_k$ with the following property: the probability that the formula 
is satisfiable (not-satisfiable) is asymptotically one (asymptotically zero) when $m/n$ is uniformly smaller (uniformly larger) than $\lambda^*$ and
$n$ diverges to infinity.
Namely, the probability of satisfiability exhibits a sharp phase transition around some constant value $\lambda^*$. The conjecture was 
positively resolved recently, first by  Carenini~\cite{GCarenini}, and  shortly after by OpenAI~\cite{k-sat-OpenAI}. 
The paper produced by OpenAI was 
one among 721 papers produced by the organization devoted to resolving some 372 open problems. 

The present paper is devoted to producing a reader-friendly exposition of the OpenAI paper. While the manuscript~\cite{k-sat-OpenAI} is fairly short
and readable
on a purely technical level, 
in a fashion typical for the works produced by  AI only, it is not optimized for the readership, contains unnecessary redundancies, applies an uncommon jargon
and obscures key ideas.
The present manuscript is devoted to fix this. Part of my motivation to write this manuscript was to understand and internalize the proof myself.

\subsection{A brief history of the conjecture}
The conjecture was formulated in Chv{\'a}tal and Reed~\cite{chvatal1992mick} who also resolved the conjecture positively for the case $k=2$.
In particular, they proved $\lambda^*_2=1$. The progress beyond $k=3$ stalled however for several decades, despite many papers devoted to the subject
and dramatic progress, particularly in narrowing the potential values of $\lambda_k^*$, assuming its validity, see for 
example~\cite{AchlioptasSorkin},\cite{AchlioptasPeres},\cite{CojaOghlanPanagiotou2016}.
Furthermore,  sharp phase transition property 
was established for the satisfiability property for all $k$ by Friedgut in his landmark paper~\cite{Friedgut}. This work produced
 necessary and sufficient conditions for the existence of sharp transition for general random graph and random hypergraph properties,
 of which satisfiability of a random k-SAT formula is a special case. 
 Unfortunately, however,
it did not resolve the conjecture itself. The import from Friedgut's main result was existence of a potentially $n$-dependent sequence $\lambda_{n,k}^*$
such that the probability of satisfiability is asymptotically one (asymptotically zero) when $m/n$ is asymptotically smaller (asymptotically larger)
than $\lambda_{n,k}^*$ and $n\to\infty$. 
In other words, while the sharpness of the transition was established, the convergence of its location was not, and as such the location
 could potentially fluctuate as a function of $n$. 

A dramatic progress towards resolving the conjecture was initiated in the field of statistical mechanics, specifically the field of spin glasses. 
Specifically, Guerra and Toninelli~\cite{GuerraTon} invented an important interpolation technique, using which they proved the existence of the ground state
energy limit for the so-called Sherrington-Kirkpatrick model. This approach was instrumental for  validating the ground state energy value
for this model, computed earlier by Parisi~\cite{parisi1979infinite,parisi1980sequence} using non-rigorous replica symmetry breaking method.
Specifically, the interpolation technique was used by Guerra~\cite{guerra2003broken} 
to prove that the value computed by Parisi  is a correct upper bound on the limiting ground state energy value. The matching lower bound
was established in a breakthrough work of Talagrand~\cite{talagrand2006parisi} and further developed in many follow up papers, see Panchenko~\cite{panchenko2013Sherrington}.

Next, closer to the random k-SAT setting, Franz and Leone~\cite{FranzLeone} adapted the interpolation technique
to the setting of diluted spin glass models, which includes random k-SAT model as a special case. A corollary of their main result is the existence of the limiting
value for the largest number of satisfiable clauses normalized by $n$ for the random k-SAT formula. The largest number of satisfiable clauses
is called the ground state
energy in the statistical physics parlance. The interpolation technique
employed in~\cite{GuerraTon} and~\cite{FranzLeone} implied superadditivity property of the expectation of the  ground state energy: the 
expected ground state
energy on the $n$-variable model is at least the sum of the expected ground state energies of the $n_1$ and $n_2$-variable models when $n=n_1+n_2$.
The existence of the limit was then a simple application of Fekete's Lemma.
Both papers relied on introducing positive temperature parameter in the respective models and then taking the temperature value to zero in the limit. 

The next development came from Bayati, Gamarnik and Tetali~\cite{BayatiGamarnikTetali} who adapted Franz and Leone's interpolation technique
directly to the zero-temperature optimization problems. The results in this paper had  implications for the satisfiability property of random
k-SAT formulas per se. Specifically,  log-superadditivity of the satisfiability property was established there, see Theorem~\ref{theorem:BGT-log-additivity} below. 
This result is at the core of the main approach developed by OpenAI in~\cite{k-sat-OpenAI}, as is explained in the main body of the present paper.
Shortly after and independently, Abbe and Montanari~\cite{abbe2014concentration} adapted Franz and Leone's technique for establishing a concentration
property for the number of satisfying assignments when conditioned on the satisfiability of the formula. They also proved in the same paper 
that a  concentration
property stronger than one they managed to obtain would imply the validity of the k-SAT conjecture. 
This is the approach successfully pursued by Carenini in~\cite{GCarenini}. A useful technical
part of~\cite{abbe2014concentration}  was a variant of the superadditivity property which was also adapted by the OpenAI solution, see 
Proposition~\ref{prop:superadditive} below. A key novelty of the Carenini's approach was obtaining a sharp bound on width of the 
Friedgut's  phase transition window. 

The next major progress towards the full resolution of the conjecture was achieved in  Ding, Sly and Sun~\cite{ding2015proof,ding2022proof}
who proved the validity of the conjecture for all sufficiently large values of $k$. Their approach was based on very refined version of the second
moment method and is technically quite involved. The upshot of this work was not only a positive resolution of the conjecture for large $k$
 but also the validation 
of the  phase transition values $\lambda_k^*$. These values were produced earlier in the statistical physics literature 
by Mertens, M{\'e}zard  and Zecchina~\cite{mertens2006threshold} also using the  replica symmetry breaking method. It is worth mentioning
that while the values produced by the replica symmetry breaking method have been confirmed in many papers such 
as~\cite{talagrand2006parisi} and~\cite{ding2015proof,ding2022proof}, the method itself remains a heuristic lacking a rigorous mathematical justification. 
The progress in~\cite{ding2015proof,ding2022proof} which took place a decade ago was the last progress until the full resolution of the conjecture
this year in~\cite{GCarenini} and~\cite{k-sat-OpenAI}.

Next we provide a brief outline of the proof. As in~\cite{BayatiGamarnikTetali} a key technical object employed by the OpenAI proof
 is the concept  of a frozen variable.
A Boolean variable $x_i$ in the  k-SAT formula is defined to be frozen when the formula is satisfiable, but every satisfying solution fixes it as the same
Boolean value, for example always at  $0$ or always at  $1$. The first technical step is then an upper bound on the number of frozen variables. 
As shown in Proposition~\ref{prop:average-density-frozen}, it is $O(n^{1\over k})$ in expectation. This bound on the expected number of frozen
variables is  used next to obtain a bound on the variance of the hitting time $H_n$ associated with random formulas. The hitting time $H_n$
is a random variable defined by considering the random k-SAT formula as a discrete time process
indexed by $m$, and constructed by adding one clause at a time. $H_n$ is then defined as largest value of $m$ for which the formula remains
satisfiable, see Subsection~\ref{subsection:preliminaries}. A bound $O(n^{1\over k})$ on the expected number of frozen variables
is  used for obtaining an upper bound $O(n^{1+{2\over k}}))$ on the variance of $H_n$, as stated in Proposition~\ref{prop:variance}.
The well-known Efron-Stein-Steele inequality Theorem~\ref{theorem:efron-stein-steele} is used towards the proof of this bound.
A direct application of Chebyshev's inequality also implies concentration property of the hitting time $H_n$, as stated in (\ref{eq:Chebyshev-tail}).
The concentration of the hitting time property allows one to relate the expected value $\E{H_n}$ of the hitting time to its  median value, namely the largest 
value of $m$ for which the likelihood of satisfiability is at least $1/2$. This is the subject of Lemma~\ref{lemma:EHn-close-lambda}.
The log-superadditivity property established in~\cite{BayatiGamarnikTetali} is then used for establishing superadditivity of the median times
as articulated in Corollary~\ref{coro:superadditive}. The final step in the proof is adaptation of the variant of the superadditivity property
developed in~\cite{abbe2014concentration}.  
This is the subject of Proposition~\ref{prop:superadditive} and this completes the description of the proof. Remarkably, as in the proof
by Carenini in~\cite{GCarenini}, the approach does not assume the validity of the Friedgut's sharp transition property and instead obtains
it as an immediate corollary. One should stress, however, that the scope of the applicability of Friedgut's theorem is much broader
than the random k-SAT formula setting. 

The remainder of the paper is organized as follows. The model description and the main theorem are stated in the next section.
The proof of the main theorem is the subject of Section~\ref{section:main-result-proof}, where the concepts of frozen variables and hitting times
are introduced in the beginning. 
Section~\ref{section:Open-problems} presents some remaining open questions. 

\section{Definitions and  the main result}
Consider $n$ Boolean variables $x_i, i\in [n]$ taking values $0$ and $1$. We fix $k$ which is clause size.
A formula $\phi$ is a conjunction of some $m$ $k$-clauses: $\phi=\wedge_{1\le j\le m} C_j$ where
each clause $C_j$ is disjunction of $k$ of the $n$ variables $x_1,\ldots,x_n$ or their negations $\bar x_1,\ldots,\bar x_n$. 
Thus for each $j$ there exists a subset $I_j^+\subset [n]$ (corresponding to variables without negation),
and a subset  $I_j^-\subset [n]$ (corresponding to negated variables), such that $| I_j^+|+| I_j^- |=k$,  $I_j^+\cap  I_j^- =\emptyset$
and
\begin{align*}
C_j=(\vee_{i\in I_j^+} x_i) \vee (\vee_{i\in I_j^-} \bar x_i)
\end{align*}
We call $I_j=I_j^+\cup I_j^-\subset [n]$ the support of the clause $C_j$. A formula $\phi$ is satisfiable if there exists
a setting of the variables $x_i$ which satisfies the formula. Namely, if there exists a vector $x=(x_1,\ldots,x_n)\in \{0,1\}^n$
such that all clauses $C_j$ evaluate to $1$ at this $x$. Such a setting, when it exists, is called a satisfying assignment or a satisfying solution. 
We note that for each clause there exists only a single assignment of variables
to $k$ values which do not satisfy the clause. When the formula is satisfiable we write $\phi\in {\rm SAT}$. When it is not satisfiable
we write $\phi\notin {\rm SAT}$.

A random formula with $n$ variables and $m$ clauses,  denoted by $\Phi(n,m)$, is obtained by generating the supports
$I_j$ uniformly at random from all ${n\choose k}$ choices, and letting  $I_j^+$ (and thus $I_j^-$) be chosen uniformly at random from $2^k$
choices as well, also independently from everything else. We denote by 
 $p(n,m)$ the probability that $\Phi(n,m)$ is sat: $p(n,m)=\pr(\Phi(n,m)\in {\rm SAT})$. The Satisfiability Conjecture asserted existence
 of a fixed threshold $\lambda^*>0$ such that $p(n,m)$ approaches $1$ when $m/n$ stays uniformly smaller than $\lambda^*$,
 and approaches $0$ when this ratio stays uniformly larger than $\lambda^*$, as $n$ diverges to infinity. 
 The OpenAI paper~\cite{k-sat-OpenAI} established precisely
 this which we now state as a formal theorem. 
 
 \begin{theorem}[Theorem 1.1, \cite{k-sat-OpenAI}]\label{theorem:main-result}
 For every $k\ge 3$ there exists $\lambda^*_k>0$ such that for every $0\le \lambda<\lambda^*_k$, $\lim_{n\to\infty} p(n,\lceil \lambda n\rceil)=1$
 and for every $\lambda>\lambda^*_k$, $\lim_{n\to\infty} p(n,\lceil \lambda n\rceil)=0$.
 \end{theorem}
 As stated in the introduction, the sharp transition result of Friedgut~\cite{Friedgut} will not be used in the proof of Theorem~\ref{theorem:main-result}.
Conversely, it is an immediate implication of Theorem~\ref{theorem:main-result}.

 \section{Proof of Theorem~\ref{theorem:main-result}}\label{section:main-result-proof}
 \subsection{Preliminaries}\label{subsection:preliminaries}
We begin by reducing the problem to the case of clauses generated with replacement. 
Specifically,  consider a modified model of the  random formula $\Phi(n,m)$ where replacement of variables is now allowed. 
Here a random
 clause $C$ is formed by selecting $k$
variables from  the list $x_1,\ldots,x_n, \bar x_1,\ldots,\bar x_n$ independently and uniformly at random \emph{with replacement}. 
The modified formula, denoted by $\Phi_{\rm wr}(n,m)$ (where ${\rm wr}$ stands for ''with replacement''), is then obtained by generating $m$ clauses with replacement
independently, as in the original version. Let $p_{\rm wr}(n,m)=\pr(\Phi_{\rm wr}(n,m)  \in {\rm SAT})$. We claim that it suffices to establish 
Theorem~\ref{theorem:main-result} with  $p_{\rm wr}(n,m)$ replacing $p(n,m)$. Specifically, suppose for every $k\ge 3$ there exists $\lambda_{k, {\rm wr}}^*$
such that for every $\lambda<\lambda^*_{k, {\rm wr}}$, $\lim_{n\to\infty} p_{\rm wr}(n,\lceil \lambda n\rceil)=1$
 and for every $\lambda>\lambda^*_{k, {\rm wr}}$, $\lim_{n\to\infty} p_{\rm wr}(n,\lceil \lambda n\rceil)=0$. Based on this assumption
 we now show the validity of Theorem~\ref{theorem:main-result} for the same limiting threshold $\lambda_k^*=\lambda_{k,{\rm wr}}^*$.
 The validity of this claim follows from a standard fact from the theory of random graphs that the collisions stemming from possible replacements
occur with probability bounded away from zero across the entire formula. Specifically, let $\mathcal{NC}$ (No-Collision) be the event that no single  clause
out of the $m$ clauses of  the random formula $\Phi_{\rm wr}(n,m)$ contains repeated variables (with or without negation). For example, 
no clause can contain a variable $x_i$ twice or variables $x_i$ and $\bar x_i$ simultaneously. It is straightforward to check that for every constant $C>0$
there exist $a>0$ such that $\liminf_n\min_{m\le Cn}\pr(\mathcal{NC}) \ge a$. 
Furthermore, conditioned on the event $ \mathcal{NC}$,
the formula $\Phi_{\rm wr}(n,m)$ is distributed as $\Phi(n,m)$. 

Now, fix any $\lambda<\lambda_{k,{\rm wr}}^*$. Then
\begin{align*}
p(n,\lceil \lambda n\rceil) & = \pr(\Phi(n,\lceil \lambda n\rceil) \in {\rm SAT}) \\
& =  \pr(\Phi_{\rm wr}(n,\lceil \lambda n\rceil) \in {\rm SAT} | \mathcal{NC})  \\
& = {  \pr(\Phi_{\rm wr}(n,\lceil \lambda n\rceil) \in {\rm SAT}, \mathcal{NC}) \over \pr(\mathcal{NC})} \\
& = {  \pr(\mathcal{NC})- \pr(\Phi_{\rm wr}(n,\lceil \lambda n\rceil) \notin {\rm SAT}, \mathcal{NC}) \over \pr(\mathcal{NC})} .
\end{align*}
We have 
\begin{align*}
\pr(\Phi_{\rm wr}(n,\lceil \lambda n\rceil) \notin {\rm SAT}, \mathcal{NC})\le \pr(\Phi_{\rm wr}(n,\lceil \lambda n\rceil) \notin {\rm SAT}) 
=1-p_{\rm wr}(n,\lceil \lambda n\rceil).
\end{align*}
By our assumption, the right-hand side converges to zero as $n\to\infty$. Since $\pr(\mathcal{NC})\ge a>0$, we conclude
$\lim_n p(n,\lceil \lambda n\rceil) =1$. 

Similarly, assume $\lambda>\lambda_{k,{\rm wr}}^*$.  Then 
\begin{align*}
p(n,\lceil \lambda n\rceil) & = \pr(\Phi(n,\lceil \lambda n\rceil) \in {\rm SAT}) \\
& =  \pr(\Phi_{\rm wr}(n,\lceil \lambda n\rceil) \in {\rm SAT} | \mathcal{NC})  \\
& = {  \pr(\Phi_{\rm wr}(n,\lceil \lambda n\rceil) \in {\rm SAT}, \mathcal{NC}) \over \pr(\mathcal{NC})} \\ 
&\le 
{  \pr(\Phi_{\rm wr}(n,\lceil \lambda n\rceil) \in {\rm SAT}) \over \pr(\mathcal{NC})} \\
&\le  {p_{\rm wr}(n,\lceil \lambda n\rceil) \over a}
\end{align*}
By our assumption $\lim_n p_{\rm wr}(n,\lceil \lambda n\rceil)=0$. We conclude $\lim_{n\to\infty}p(n,\lceil \lambda n\rceil)=0$. 
This completes the proof of the claim.

The remainder of the paper is thus devoted to establishing Theorem~\ref{theorem:main-result} for $p_{\rm wr}(n,m)$ replacing $p(n,m)$.
For simplicity,  we denote $\Phi_{\rm wr}(n,m)$ and $p_{\rm wr}(n,m)$  by $\Phi(n,m)$ and $p(n,m)$,
respectively, throughout the rest of the paper.
 
 Next we introduce the concepts of  hitting times and frozen variables. We begin by defining hitting times.
 Given $n$ and $m$ let $Z$ denote the (random) number of  solutions satisfying the formula $\Phi(n,m)$. A straightforward computation shows
 $\E{Z}=2^n(1-1/2^k)^m$. By Markov inequality $p(n,m)\le \E{Z}$. When $m\ge 2^{k}n$, we obtain 
 \begin{align*}
 p(n,m)\le 2^n(1-1/2^k)^m\le (2e^{-1})^n<1/2
 \end{align*}
 for all sufficiently large $n$. Define $m(n)$ as the  largest $m$ such that $p(n,m)\ge 1/2$. Then  $p(n,m(n)+1)<1/2$.  By the above
 $m(n)<2^{k}n$ for all $n$.
Let  $\lambda_n=m(n)/n$. Then for all $n$, $\lambda_n<2^{k}$. 
In particular,  $\lambda_n=O(1)$. It is well known that for all $\epsilon>0$, $\liminf_n \lambda_n\ge 1-\epsilon$ and thus for some $n_0$ and all $n\ge n_0$
\begin{align}\label{eq:lambda-order-1}
1/2\le \lambda_n\le 2^{k}.
\end{align}
See for example~\cite{ChvatalReed} which shows that for $k=2$, $\lambda_n\to 1$ with high concentration around $\lambda_n$
for the version without replacement. By the equivalence property above the same applies to the model with replacement.
Finally, we use $\lambda_{n,k}\le \lambda_{n,k+1}$ to assert (\ref{eq:lambda-order-1}).

Let 
$H_n$ be the hitting time defined as follows. Consider the random formula $\Phi(n,2^{k+1}n)$ and the associated stochastic process
$(\Phi(n,m), 0\le m\le 2^{k+1}n)$.
Define
\begin{align*}
H_n=\max\{m\le 2^{k+1}n: \Phi(n,m)\in {\rm SAT}\}.
\end{align*}
Namely $H_n$ is the largest index of a  clause  which  maintains the satisfiability property, capped at $2^{k+1}n$.
Then for all $m\le 2^{k+1}n$ we have $p(n,m)=\pr(H_n\ge m)$. Thus $m(n)$ is the median of the random variable $H_n$ and $\lambda_n$
is median normalized by $n$.
 
Next we define frozen variables. Given a formula $\phi$ defined on variables $x_1,\ldots,x_n$, a variable $x_i$ is called
frozen if it takes the same value for all satisfying assignments. Namely $x_i, i\in [n]$ is  frozen if either $x_i=0$ for all satisfying solutions
$(x_1,\ldots,x_n)\in \{0,1\}^n$, or $x_i=1$ for all satisfying solutions $(x_1,\ldots,x_n)\in \{0,1\}^n$.
If the formula is not satisfiable, (and thus the number of satisfying assignments is empty), the set of frozen variables is empty. 
Given a satisfiable formula $\phi$ and a clause $C$ note that the modified formula $\phi \wedge C$
is not satisfiable if and only if every variable in the support of $C$ is frozen and its parity agrees with the value of the frozen variable. 
For example, if $C=(x_3\vee \bar x_7 \vee \bar x_8)$ then $\phi\wedge C$ is not satisfiable if and only if variables $x_3,x_7,x_8$ 
are frozen to values $0,1,1$ respectively. 

Let $0\le f(\phi)\le n$ denote the number of frozen variables in formula $\phi$. We adopt the convention that $f(\phi)=0$
when $\phi$ is not satisfiable.
Note that for every 
formula $\phi$, $(f(\phi)/(2n))^k$ is precisely the probability that the modified formula $\phi \wedge C$ is not satisfiable 
when $C$ is a $k$-clause chosen uniformly at random. We write $q_k(\phi)=(f(\phi)/(2n))^k$ for short.

\ignore{
When $f(\phi)\ge k$ define 
\begin{align}
q_k(\phi)={\prod_{0\le r\le k-1} (f(\phi)-r) \over 2^k\prod_{0\le r\le k-1} (n-r)} \le \left( {f(\phi) \over 2n}\right)^k. \label{eq:f-in-q}
\end{align}
}

\ignore{
We set $q_k(\phi)=0$ when $f(\phi)\le k-1$. This is consistent with the fact that adding any single $k$-clause to the formula maintains the
satisfiability property since the 
support of $C$ consists of $k$ distinct elements.
}

\subsection{Controlling frozen variables}
The goal of this subsection is to obtain an upper bound on the number of frozen variables in a random formula $\Phi(n,m)$. In particular,
we will establish in Proposition~\ref{prop:average-density-frozen} below that the expected  number of the frozen variables in formula $\Phi(n,m)$
is  order $n^{1\over k}$ \emph{on average} when we average over the number of clauses $m$ in the range $[0,2^{k+1}n]$.

We begin with a preliminary lemma. It shows that, modulo order $n^{-{k-1\over k}}$ correction term,  
for every satisfiable formula $\phi$, appending this formula with a $k-1$ clause chosen
uniformly at  random leads to a higher chance of maintaining satisfiability vs appending it with order $n^{1\over k}$ many $k$-clauses also 
chosen uniformly at random and independently. 
\begin{lemma}\label{lemma:k-1-to-k-clauses}
Let $\phi$ be any fixed  formula on $n$ variables. Let $C_0$ be a $k-1$-clause chosen uniformly at random and 
let $C_1,\ldots, C_R$ be $k$-clauses chosen independently uniformly at random. Suppose $R\ge 2n^{1\over k}$.
There exists sufficiently large $n_0$ such that for all $n\ge n_0$
\begin{align}
\pr(\phi \wedge C_0\in {\rm SAT})\ge \pr\left(\phi \wedge  \left(\wedge_{1\le j\le R} C_j\right) \in {\rm SAT}\right)-n^{-{k-1\over k}}. \label{eq:k-1-to-k-clauses}
\end{align}
\end{lemma}

\begin{proof}
The claim holds trivially when $\phi$ is not satisfiable. Thus assume $\phi$ is satisfiable.
Suppose first $f(\phi)/(2n)\le n^{-{1\over k}}$. Then $\pr(\phi \wedge C_0\in {\rm SAT})=1-(f(\phi)/(2n))^{k-1}\ge 1-n^{-{k-1\over k}}$, and the claim holds trivially.

Thus suppose $f(\phi)/(2n)\ge n^{-{1\over k}}$. Then
\begin{align*}
R(f(\phi)/(2n))^k & \ge 2n^{1\over k} (f(\phi)/(2n))^k \\
&= 2n^{1\over k} (f(\phi)/(2n)) (f(\phi)/(2n))^{k-1} \\
&\ge 2 (f(\phi)/(2n))^{k-1}.
\end{align*}
Then
\begin{align*}
\pr\left(\phi \wedge  \left(\wedge_{1\le j\le R} C_j\right) \in {\rm SAT}\right) 
& \le 
\pr\left(\cap_{1\le j\le R} \left( \phi \wedge C_j \in {\rm SAT}\right)\right) \\
&= \left(1-\left({f(\phi)\over 2n}\right)^k\right)^R \\
&\le e^{-R(f(\phi)/(2n))^k} \\
&\le e^{-2(f(\phi)/(2n))^{k-1}} \\
&\le 1-(f(\phi)/(2n))^{k-1} \\
&=\pr(\phi \wedge C_0\in {\rm SAT}),
\end{align*}
where we use $(f(\phi)/(2n))^k \le 1/2^{k-1}\le 1/2$ and $e^{-2x}\le 1-x$ for $x\in [0,1/2]$.
\end{proof}

The following bound shows that 
for random formulas on $n$ variables the number of frozen variables is at most order $n^{1\over k}$ ''on average'' when we average
over $m$ in the range $[0,2^{k+1}n]$.

\begin{prop}\label{prop:average-density-frozen}
For every $k$ there exists a constant $c_k$ and $n_0$ such that for all $n\ge n_0$
\begin{align*}
\sum_{0\le m\le 2^{k+1} n}{ \E{f(\Phi(n,m))} \over n} \le c_k n^{1\over k}.
\end{align*}  
\end{prop}

\begin{proof}
Fix $m\le 2^{k+1}n$.
By symmetry $\E{f(\Phi(n,m))}/n$ equals the probability that $\Phi(n,m)$ is satisfiable and the  variable $x_n$ is frozen, since when
$\Phi(n,m)$ is not satisfiable, $f(\Phi(n,m))=0$ by assumption.
Let $\Phi_{-n}$ be the part of this formula spanned by the first $n-1$ variables $x_1,\ldots,x_{n-1}$.
Consider the clauses $C_1,\ldots,C_D$ containing the variable $x_n$ (where $D$ is random), so that $\Phi(n,m)=\Phi_{-n}\wedge (\wedge_{1\le j\le D} C_j)$.
Consider the clauses $\hat C_j, 1\le j\le D$ obtained from $C_j, 1\le j\le D$ by deleting variable $x_n$. 
Each clause has at  most $k-1$ variables. Specifically, the clause size $\hat C_j$ is strictly less than $k-1$ if variable $x_n$ appears in $C_j$ more
than once.  Let $\mathcal{C}_{k-1}$ be the event that all  clauses $\hat C_j, 1\le j\le D$ have size $k-1$. 
We can find constant $c>0$ such that $\pr(\mathcal{C}_{k-1})\ge 1-c/n$  uniformly in $m\le 2^{k+1}n$.

Next, we  make the following observation: if $x_n$ is a frozen variable in $\Phi(n,m)$  then $\Phi_{-n}$ is satisfiable (since $\Phi(n,m)$ is), but
the formula obtained by appending $\hat C_j, 1\le j\le D$ is not satisfiable. That is $\Phi_{-n}\in {\rm SAT}$, but 
$\hat \Phi_{-n}\triangleq \Phi_{-n}\wedge (\wedge_{1\le j\le D} \hat C_j)\notin {\rm SAT}$. Indeed, otherwise a satisfying assignment of $\hat \Phi(n,m)$ can be extended
to a satisfying assignment of $\Phi(n,m)$ with either choice $x_n=0$ or $x_n=1$, contradicting the assumption that $x_n$ is frozen. Thus
\begin{align}\label{eq:frozen-upper-bound}
{\E{f(\Phi(n,m))}\over n}=\pr(x_n~{\rm frozen})\le \pr(\Phi_{-n}\in {\rm SAT}, \hat\Phi_{-n}\notin {\rm SAT}).
\end{align}
By earlier observation, 
\begin{align}
\pr(\Phi_{-n}\in {\rm SAT}, \hat\Phi_{-n}\notin {\rm SAT})\le \pr(\Phi_{-n}\in {\rm SAT}, \hat\Phi_{-n}\notin {\rm SAT} | \mathcal{C}_{k-1})+O(1/n).
\label{eq:condition-on-k-1}
\end{align}

We now construct a sequence of modified formulas $\hat \Phi_j, 0\le j\le D$ starting 
with $\hat \Phi_0=\hat \Phi_{-n}$ as follows. For $j=1,2,\ldots,D$, the formula $\hat\Phi_j$ is obtained
from $\hat\Phi_{j-1}$ by removing clause $\hat C_j$ and replacing it with $R\triangleq \lceil 2n^{1\over k} \rceil$  $k$-clauses $C_{j,1},\ldots,C_{j,R}$
on variables $x_1,\ldots,x_{n-1}$ and their negations.
These clauses
are generated uniformly at random independently and also independently from the all the other clauses. 

Fix $0\le d\le m$. 
We note that conditionally on $D=d$,  $\hat \Phi_D$ is distributed as  $\Phi(n-1,m-d+Rd)$. 
Applying Lemma~\ref{lemma:k-1-to-k-clauses} we have for some constant $n_0$ for all $n\ge n_0$ and for all $j=0,\ldots, d-1$
\begin{align*}
\pr(\hat\Phi_j \in {\rm SAT} | D=d, \mathcal{C}_{k-1})\ge \pr\left( \hat\Phi_{j+1}\in {\rm SAT}  | D=d, \mathcal{C}_{k-1} \right)-(n-1)^{-{k-1\over k}}. 
\end{align*}
Summing over $j$ we obtain that conditional on $D=d$
\begin{align*}
\pr(\hat\Phi_{-n} \in {\rm SAT} | D=d, \mathcal{C}_{k-1})\ge p(n-1,m-d+Rd)-d(n-1)^{-{k-1\over k}}. 
\end{align*}
On the other hand 
\begin{align*}
\pr(\hat\Phi_{-n} \in {\rm SAT} | D=d, \mathcal{C}_{k-1})=\pr(\Phi_{-n}\in {\rm SAT}|D=d, \mathcal{C}_{k-1})
-\pr(\Phi_{-n}\in {\rm SAT}, \hat\Phi_{-n}\notin {\rm SAT} |D=d, \mathcal{C}_{k-1}).
\end{align*}
We obtain
\begin{align*}
&\pr(\Phi_{-n}\in {\rm SAT}, \hat\Phi_{-n}\notin {\rm SAT} |D=d,  \mathcal{C}_{k-1}) \\
&\le \pr(\Phi_{-n} \in {\rm SAT} | D=d, \mathcal{C}_{k-1})-p(n-1,m-d+Rd)+d(n-1)^{-{k-1\over k}}.
\end{align*}
We note 
\begin{align*}
\pr(\Phi_{-n} \in {\rm SAT} | D=d, \mathcal{C}_{k-1})=p(n-1,m-d). 
\end{align*}
We conclude
\begin{align*}
\pr(\Phi_{-n}\in {\rm SAT}, \hat\Phi_{-n}\notin {\rm SAT} |D=d, \mathcal{C}_{k-1})\le p(n-1,m-d)-p(n-1,m-d+Rd)+d(n-1)^{-{k-1\over k}}.
\end{align*}
Summing over $m=0,1,\ldots,2^{k+1}n$ we obtain
\begin{align*}
\sum_{d\le m\le 2^{k+1}n}\pr(\Phi_{-n}\in {\rm SAT}, \hat\Phi_{-n}\notin {\rm SAT} |D=d, \mathcal{C}_{k-1})
&\le \sum_{0\le m\le Rd-1}p(n-1,m)+2^{k+1}n d(n-1)^{-{k-1\over k}} \\
&\le 
Rd+2^{k+1}dn^{{1\over k}} \\
&\le 
c_{k,1}dn^{1\over k},
\end{align*}
for some $k$ dependent constant $c_{k,1}$. Multiplying throughout by likelihood that $D=d$ we obtain
\begin{align*}
\sum_{d\le m\le 2^{k+1}n}\pr(\Phi_{-n}\in {\rm SAT}, \hat\Phi_{-n}\notin {\rm SAT}, D=d | \mathcal{C}_{k-1})
&\le 
c_{k,1}n^{1\over k}d\max_{d\le m\le 2^{k+1}n}\pr_m(D=d  | \mathcal{C}_{k-1}),
\end{align*}
where subscript  in $\pr_m$ emphasizes  that the probability depends on $m$. 
For fixed $m$ the random variable $D$ has Binomial distribution with  $m=O(n)$ trials and success probability $O(1/n)$. 
It is easy to see that  $\sum_{d\ge 0}d\max_{d\le m\le 2^{k+1}n}\pr_m(D=d  | \mathcal{C}_{k-1})\le c_{k,2}<\infty$. 
 for some constant $c_{k,2}$. Summing over $d$, we obtain
\begin{align*}
\sum_{0\le m\le 2^{k+1}n}\pr(\Phi_{-n}\in {\rm SAT}, \hat\Phi_{-n}\notin {\rm SAT} | \mathcal{C}_{k-1})
&\le 
c_{k,3}n^{1\over k},
\end{align*}
and further using (\ref{eq:condition-on-k-1}) we obtain
\begin{align*}
\sum_{0\le m\le 2^{k+1}n}\pr(\Phi_{-n}\in {\rm SAT}, \hat\Phi_{-n}\notin {\rm SAT})
&\le 
c_{k,3}n^{1\over k}+c_{k,4}
\end{align*}
for some constant $c_{k,4}$. Combining with (\ref{eq:frozen-upper-bound}) and (\ref{eq:condition-on-k-1})
we complete the proof by increasing $c_k$ to absorb the additive constant $c_{k,4}$.
\end{proof}

\subsection{Concentration of the hitting times $H_n$}
The goal of this subsection is to establish a very important technical result: the hitting time $H_n$ is concentrated around its expectation.
Specifically, in Proposition~\ref{prop:variance} below we will show that the variance of $H_n$ is at most order $n^{1+{2\over k}}$. The concentration then will 
follow directly by applying Chebyshev's inequality.

Given formula $\Phi(n,2^{k+1}n)$ and index $1\le j\le 2^{k+1}n$ let $\Phi_{-j}$ denote the formula $\Phi(n,2^{k+1}n)$ with clause $C_j$ omitted. Let $H_{-j}$ be the corresponding
hitting time. To be precise, let $C_1,\ldots,C_{2^{k+1}n}$ be the clauses of $\Phi(n,2^{k+1}n)$. Denote by $\Phi_{m,-j}$ the formula $\wedge_{1\le \ell \le m, \ell\ne j} C_\ell$.
Then
\begin{align*}
H_{-j}=\max\{m\le 2^{k+1}n: \Phi_{m,-j}\in {\rm SAT}\}.
\end{align*}

We omit index $n$ for brevity. Also we write $H$ for $H_n$.  Trivially $H\le H_{-j}$ for all $j\ge 1$. 
The following bound on the variance of $H$ is based on standard inequality known as Efron-Stein and Steele's application of it 
to non-symmetric cases. We include the proof for completeness. 

\begin{theorem}[Efron-Stein-Steele bound]\label{theorem:efron-stein-steele}
The following holds for all $n$:
\begin{align*}
{\rm Var}(H)\le  \sum_{1\le j\le 2^{k+1}n}\E{(H-H_{-j})^2}.
\end{align*}
\end{theorem}

\begin{proof}
For $j=0,1\ldots,2^{k+1}n$ denote by $\mathcal{F}_j$ the $\sigma$-field generated by randomness of $C_1,\ldots,C_j$, where $\mathcal{F}_0$ is assumed to be the trivial
$\sigma$-field. Note that $H$ is measurable with respect to $\mathcal{F}_{2^{k+1}n}$, since the realizations of $C_j, 1\le j\le 2^{k+1}n$ 
uniquely determine the formula and the value
of $H$ as a result. Introduce martingale differences for $1\le j\le 2^{k+1}n$
\begin{align*}
D_j=\E{H | \mathcal{F}_j}-\E{H | \mathcal{F}_{j-1}}
\end{align*}
so that $H-\E{H}=\sum_{1\le j\le 2^{k+1}n} D_j$. A standard property  satisfied by martingale differences is the orthogonality: for all $j\ne \ell$ $\E{D_jD_\ell}=0$.
It is a simple consequence of the tower property $\E{\E{H | \mathcal{F}_j} | \mathcal{F}_\ell}=\E{H | \mathcal{F}_\ell}$ for all $j\ge \ell$. Orthogonality implies
\begin{align}
{\rm Var}(H) = \sum_{1\le j\le 2^{k+1}n} \E{D_j^2}. \label{eq:sum-Dj}
\end{align}
Next, as $H_{-j}$ does not depend on $C_j$ we have $\E{H_{-j} | \mathcal{F}_j}=\E{H_{-j} | \mathcal{F}_{j-1}}, j\ge 1$ which means
\begin{align*}
D_j=\E{H-H_{-j} | \mathcal{F}_j}-\E{H-H_{-j} | \mathcal{F}_{j-1}}.
\end{align*}
Applying the tower property we obtain
\begin{align*}
\E{D_j^2}=\E{\left(\E{H-H_{-j} | \mathcal{F}_j}\right)^2} -\E{\left(\E{H-H_{-j} | \mathcal{F}_{j-1}}\right)^2} \le \E{\left(\E{H-H_{-j} | \mathcal{F}_j}\right)^2}.
\end{align*}
Applying conditional Jensen's inequality $\E{D_j^2}\le \E{(H-H_{-j})^2}$. Combining with (\ref{eq:sum-Dj}) we complete the proof.
\end{proof}

The next lemma relates the variance of the hitting time to the expectation of the number of frozen variables. 
\begin{lemma}\label{lemma:var-frozen}
The following holds for every $j=1,2,\ldots,2^{k+1}n$
\begin{align*}
\E{(H-H_{-j})^2}\le 2\sum_{j\le m\le 2^{k+1}n} \E{\min(2^{k+1}nq_k(\Phi(n,m-1)),1)}.
\end{align*}
\end{lemma}

\begin{proof}
Fix $1\le j\le 2^{k+1}n$.
Recall our earlier notation $\Phi_{m,-j}$ for $m\ge j$. This is the formula obtained from the first $m$ clauses of the clauses $C_1,\ldots,C_{2^{k+1}n}$ of 
$\Phi(n,2^{k+1}n)$ with clause $C_j$ omitted. 
Let $\Phi_m$ be the same formula without omitting $C_j$, so that $\Phi_m$ is $\Phi(n,m)$ in distribution. Introduce the events
\begin{align*}
\mathcal{A}_m=\{ \Phi_{m,-j}\in {\rm SAT}, \Phi_m\notin {\rm SAT}\}.
\end{align*}
This is precisely the event $H<m\le H_{-j}$. Then
\begin{align*}
H_{-j}-H= | \{m: H<m\le H_{-j}\} |=\sum_{j\le m\le 2^{k+1}n} {\bf 1}(\mathcal{A}_m). 
\end{align*} 
We obtain 
\begin{align}
\E{(H_{-j}-H)^2} &= \E{ \left(\sum_{j\le m\le 2^{k+1}n} {\bf 1}(\mathcal{A}_m)\right)^2} \notag\\
&=
\sum_{j\le m\le  2^{k+1}n} \pr(\mathcal{A}_m)+
2\sum_{j\le m < \ell \le 2^{k+1}n} \pr(\mathcal{A}_m\cap \mathcal{A}_\ell )   \notag\\
&\le 
2\sum_{j\le m\le \ell \le 2^{k+1}n} \pr(\mathcal{A}_m\cap \mathcal{A}_\ell ) . \label{eq:2-sum-m-l}
\end{align} 
Next we claim that for every $\ell \ge m$
\begin{align*}
\mathcal{A}_m\cap \mathcal{A}_\ell=\mathcal{A}_m\cap \{\Phi_{\ell,-j}\in {\rm SAT}\}.
\end{align*}
Indeed, since $\mathcal{A}_\ell\subset \{\Phi_{\ell,-j}\in {\rm SAT}\}$ then 
\begin{align*}
\mathcal{A}_m\cap \mathcal{A}_\ell\subset \mathcal{A}_m\cap \{\Phi_{\ell,-j}\in {\rm SAT}\}.
\end{align*}
On the other hand $\mathcal{A}_m\subset \{\Phi_m\notin {\rm SAT}\}\subset \{\Phi_{\ell}\notin {\rm SAT}\}$, the latter following from $\ell\ge m$. 
This furnishes the other inclusion.

Consider the $\sigma$-field $\mathcal{F}_{m,-j}$ generated by clauses $C_r, 1\le r\le m, r\ne j$.  
The randomness of clauses $C_j$ as well as clauses
$C_{m+1},\ldots,C_{2^{k+1}n}$ is independent from this $\sigma$-field. Recall our notation $q_k(\phi)=(f(\phi)/(2n))^k$.
Then 
$\pr(\mathcal{A}_m | \mathcal{F}_{m,-j})=q_k(\Phi_{m,-j})$. Also $\Phi_{\ell,-j}\in {\rm SAT}$ implies $\Phi_{m,-j}\wedge C_r \in {\rm SAT}$
for all $m+1\le r\le \ell$. 
Similarly, we have  $\pr(\Phi_{m,-j}\wedge C_r \in {\rm SAT} | \mathcal{F}_{m,-j})=1-q_k(\Phi_{m,-j})$ for $r\ge m+1$.
The events $\mathcal{A}_m$ and $\Phi_{m,-j}\wedge C_r\in {\rm SAT}, r=m+1,\ldots,\ell$ are independent conditionally on $\mathcal{F}_{m,-j}$,
since they depend on randomness of disjoint clauses, namely $C_j$ and $C_r, m+1,\ldots,\ell$. 
Thus
\begin{align*}
\pr(\mathcal{A}_m\cap \mathcal{A}_\ell | \mathcal{F}_{m,-j}) 
&=
\pr(\mathcal{A}_m\cap \{\Phi_{\ell,-j}\in {\rm SAT}\}  | \mathcal{F}_{m,-j} ) \\
&\le  q_k(\Phi_{m,-j})(1-q_k(\Phi_{m,-j}))^{\ell-m}.
\end{align*}
Then 
\begin{align*}
\sum_{\ell:  m\le \ell \le 2^{k+1}n}\pr(\mathcal{A}_m\cap \mathcal{A}_\ell | \mathcal{F}_{m,-j}) 
&\le  \sum_{\ell:  m\le \ell \le 2^{k+1}n}q_k(\Phi_{m,-j})(1-q_k(\Phi_{m,-j}))^{\ell-m} \\
&\le  \min( 2^{k+1}n q_k(\Phi_{m,-j}), 1).
\end{align*}
Taking the outer expectation with respect to $\mathcal{F}_{m,-j}$ we obtain
\begin{align*}
\sum_{\ell:  m\le \ell \le 2^{k+1}n}\pr(\mathcal{A}_m\cap \mathcal{A}_\ell ) 
\le \E{\min( 2^{k+1}n q_k(\Phi_{m,-j}), 1)}.
\end{align*}
Summing over $m=j,j+1,\ldots,2^{k+1}n$ and using (\ref{eq:2-sum-m-l}) we conclude the proof.
\end{proof}

We are ready to state and prove the concentration of hitting times property. 
\begin{proposition}\label{prop:variance}
For every $k$ there exists large enough constant $C_k$ such that ${\rm Var}(H_n)\le C_kn^{1+{2\over k}}$ for all $n$. As a result
for $\delta\triangleq (k-2)/(4k)$ and $\eta=(k-2)/(2k)$
\begin{align}
\pr\left( | H_n-\E{H_n} |\ge n^{1-\delta}\right)\le C_k n^{-\eta} \label{eq:Chebyshev-tail}.
\end{align}
\end{proposition}
While $\delta$ and $\eta$ were simple to express in terms of $k$, their detailed values are not of importance to the remainder of the proof. 
As we will see, it suffices to have $0<\delta<1$ and $\eta>0$. The only adjustment is that the exponent $\delta$ appearing
in Corollary~\ref{coro:superadditive} in the next subsection would need to be reduced to a smaller positive value. 

\begin{proof}
The tail bound (\ref{eq:Chebyshev-tail}) follows from the bound on variance ${\rm Var}(H_n)$ immediately by Chebyshev's inequality. Thus we focus 
on bounding the variance ${\rm Var}(H_n)$. Using $\min(z,1)\le z^{1\over k}$ for all $z\ge 0$, we have
\begin{align*}
\min\left(2^{k+1}n q_k(\Phi(n,m)),1\right) \le 2^{1+1/k}n^{1\over k} (q_k(\Phi(n,m)))^{1\over k} =
2^{1+{1\over k}}n^{1\over k} (f(\Phi(n,m))/n).
\end{align*}
Applying Lemma~\ref{lemma:var-frozen} we obtain
\begin{align*}
\E{(H-H_{-j})^2}\le 4n^{1\over k}\sum_{j\le m\le 2^{k+1}n} n^{-1} \E{f(\Phi(n,m-1))}.
\end{align*}
Applying Proposition~\ref{prop:average-density-frozen} we obtain
\begin{align*}
\E{(H-H_{-j})^2}\le 4c_kn^{2\over k}.
\end{align*}
Finally, using Efron-Stein-Steele bound Theorem~\ref{theorem:efron-stein-steele} we obtain

\begin{align*}
{\rm Var}(H)\le  4c_k2^{k+1}n^{1+{2\over k}}
\end{align*}
This completes the proof by setting $C_k=4c_k2^{k+1}$.
\end{proof}

\subsection{From concentration of  $H_n$ to superadditivity of the normalized medians}
The main technical result for this section is establishing the superadditivity type property for the normalized
median times $\lambda_n$ which is subject of Corollary~\ref{coro:superadditive},
as well as continuity of the median times in $n$. Specifically, in Lemma~\ref{lemma:from-n-to-n+1} we establish that $\lambda_{n+1}$ is close to $\lambda_n$.
We begin by showing the proximity of mean and median of the hitting times.

\begin{lemma}\label{lemma:EHn-close-lambda}
Let $\delta=(k-2)/(4k)$ be as in Proposition~\ref{prop:variance}.
There exists $n_0$ such that for all $n\ge n_0$
\begin{align*}
|\E{H_n}-\lambda_n n|\le 2n^{1-\delta}.
\end{align*}
\end{lemma}
As we will see in the proof, the choice of multiplier $2$ is somewhat arbitrary and any constant strictly larger than $1$ suffices.
\begin{proof}
Suppose first $\E{H_n}\ge \lambda_n n+2n^{1-\delta}$. Then by Proposition~\ref{prop:variance} $\pr(H_n\le \lambda_n n+n^{1-\delta})\le n^{-\eta}$. 
On the other hand 
\begin{align*}
{1\over 2}>p(n,\lambda_n n+1)=\pr(H_n\ge \lambda_n n+1),
\end{align*}
implying 
\begin{align*}
\pr(H_n\le \lambda_n n+n^{1-\delta})\ge \pr(H_n\le \lambda_n n+1)> 1/2>n^{-\eta},
\end{align*} 
for large enough $n$ which is a contradiction. We conclude $\E{H_n}\le \lambda_n n+2n^{1-\delta}$.

Conversely, suppose $\E{H_n}\le \lambda_n n-2n^{1-\delta}$. Then $\pr(H_n\ge \lambda_n n-n^{1-\delta})\le n^{-\eta}$. 
On the other hand 
\begin{align*}
{1\over 2}\le p(n,\lambda_n n)=\pr(H_n\ge \lambda_n n)\le \pr(H_n\ge \lambda_n n-n^{1-\delta}),
\end{align*}
which is again a contradiction when $n^{-\eta}<1/2$. We conclude $\E{H_n}\ge \lambda_n n-2n^{1-\delta}$.
\end{proof}

Next we import the following superadditivity result from~\cite{BayatiGamarnikTetali} based on the  combinatorial interpolation technique. 
We reproduce the proof for completeness.
This  result is at the core of establishing the satisfiability limit as stated in our main result, Theorem~\ref{theorem:main-result}.
\begin{theorem}\label{theorem:BGT-log-additivity}
Given $n, n_1, n_2$ such that $n=n_1+n_2$ and given $m$, let random variables $M, M_1, M_2$ have
 Poisson distributions with means $m, (n_1/n)m, (n_2/n)m$, respectively. 
Then  
\begin{align*}
\E{p(n,M)}\ge \E{p(n_1,M_1)}\E{p(n_2,M_2)} ,
\end{align*}
where  the expectations of $p(n,M)$ and $p(n_j,M_j), j=1,2$  are with respect to the randomness of $M,M_j, j=1,2$.
\end{theorem}

\begin{proof}
Consider the random formula $\Phi(n,M)$. 
On the event $M>0$ we consider the following process indexed by time $t=1,\ldots,M$. We order $M$ clauses uniformly at random and independently from
any other randomness of the formula. For each $t$ we delete the $t$-th clause $C_t$. We then generate a new $k$-clause $\hat C_t$ as follows. 
With probability $n_1/n$ this clause is generated uniformly at random with replacement from the set of variables $x_1,\ldots,x_{n_1}$ and their negations. 
With the remaining probability
$n_2/n$ the clause is generated uniformly at random with replacement from  the set of variables $x_{n_1+1},\ldots,x_{n}$ and their negations. 
If $M=0$ we do nothing.
Let $\Phi_t, t=0,1,\ldots,M$ be the random formula obtained at time $t$, after clause $C_t$ is deleted and $\hat C_t$ is created.
In particular,  $\Phi_0=\Phi(n,M)$.
Also by  the Poisson splitting theorem, the  formula $\Phi_M$ at the end of the process is a disjoint union of independent formulas $\Phi_j, j=1,2$ on the two sets of nodes
$1,\ldots,n_1$ and $n_1+1,\ldots,n$, respectively,
where $\Phi_j$ is distributed as $\Phi(n_j,M_j)$ for $j=1,2$.

Fix $T>0$ and condition on $M=T$. We claim that for every  $1\le t\le T$ 
\begin{align}
\pr(\Phi_{t-1}\in {\rm SAT} | M=T)\ge \pr(\Phi_t\in {\rm SAT} | M=T), \label{eq:monotonicity-Phi-t}
\end{align}
from which the assertion of the theorem follows.

Let $\Delta\Phi_t$ be the random formula obtained after deleting clause $C_t$ but before adding the clause $\hat C_t$.
Then $\Phi_{t-1}=\Delta\Phi_t \wedge C_t$ and $\Phi_{t}=\Delta\Phi_t \wedge \hat C_t$.
We condition on the realization $\phi$ of $\Delta\Phi_t$. Note that $C_t$ and $\hat C_t$ are independent from $\phi$. Assume first $\phi$ is  satisfiable formula. Then
$\pr(\phi\wedge C_t \notin {\rm SAT})=\left({f(\phi) \over 2 n}\right)^k.$ For $j=1,2$, let $f_j(\phi)$ be the number of frozen variables of the formula $\phi$ in subsets $1,\ldots,n_1$ and
$n_1+1,\ldots,n$, respectively. Then 
\begin{align}\label{eq:convexity-1}
\pr(\phi\wedge \hat C_t \notin {\rm SAT})={n_1\over n}\left({f_1(\phi) \over 2 n_1}\right)^k+{n_2\over n}\left({f_2(\phi) \over 2 n_2}\right)^k.
\end{align}
Using the convexity of the function $x^k$ when restricted to the domain $x\ge 0$, the right-hand side is at least
\begin{align}\label{eq:convexity-2}
\left({n_1\over n}{f_1(\phi) \over 2 n_1}+{n_2\over n}{f_2(\phi) \over 2 n_2}\right)^k = \left({ f_1(\phi)+f_2(\phi) \over 2n }\right)^k=\pr(\phi\wedge C_t \notin {\rm SAT}). 
\end{align}
We conclude $\pr(\phi\wedge \hat C_t \notin {\rm SAT})\ge \pr(\phi\wedge C_t \notin {\rm SAT})$. 

On the other hand, when  formula $\phi$ is not satisfiable, we have $\pr(\phi\wedge \hat C_t \in {\rm SAT})=\pr(\phi\wedge C_t \in {\rm SAT})=0$. 
We conclude that for every realization of $\phi$, $\pr(\phi\wedge  C_t \in {\rm SAT})\ge \pr(\phi\wedge \hat C_t \in {\rm SAT})$. 
Thus (\ref{eq:monotonicity-Phi-t}) is verified.
Removing the conditioning on $\phi$ and $M=T$ we complete the proof. 
\end{proof}

The key implication of Theorem~\ref{theorem:BGT-log-additivity} is the following property of the normalized median times. 
\begin{coro}\label{coro:superadditive}
There exists $n_0$ such that for all  $n_1,n_2\ge n_0$ 
\begin{align*}
\lambda_{n_1+n_2}\ge \min(\lambda_{n_1},\lambda_{n_2})-5 (n_1+n_2)^{-\delta}.
\end{align*}
\end{coro}

\begin{proof}
Fix $n_1, n_2>1$ and let $n=n_1+n_2$.
Let $\lambda=\min_{j=1,2} \lambda_{n_j}$. Set $r=\lambda n$. Then $r_j\triangleq (n_j/n)r\le \lambda_{n_j}n_j$.
(We note that  $r, r_j$ are not necessarily integers, unlike $\lambda_nn, \lambda_{n_j}n_j, j=1,2$,) 
Let $M,M_j, j=1,2$ be Poisson random variables with means $r,r_j, j=1,2$ respectively. We use the following
 universal lower bound for a Poisson random variable $X$ with mean $\mu$: $\pr(X\le \mu)\ge 1/e$.
This implies for $j=1,2$
\begin{align*}
\E{p(n_j,M_j)} & \ge \pr(\Phi(n_j, M_j)\in {\rm SAT} | M_j\le r_j)\pr( M_j\le r_j) \\
&\ge 
p(n_j,\lambda_{n_j}n_j)(1/e)\\
&\ge
1/(2e).
\end{align*}
Applying  Theorem~\ref{theorem:BGT-log-additivity}, we obtain $\E{p(n,M)}\ge 1/(2e)^2$. 
Fix $\gamma\in (1/2,1-\delta)$. This is possible since $\delta=(k-2)/(4k)<1/2$ for all $k$.
Since ${\rm Var}(M)=\lambda n$, by Chebyshev's inequality we have 
\begin{align*}
\pr(M-r < -n^{\gamma})\le {\lambda n\over n^{2\gamma}}={\lambda \over n^{2\gamma-1}}.
\end{align*}
Thus
\begin{align*}
\pr(\Phi(n,M)\in {\rm SAT} | M\ge r-n^{\gamma})\pr(M\ge r-n^{\gamma})\ge 1/(2e)^2-{\lambda \over n^{2\gamma-1}},
\end{align*}
namely
\begin{align*}
p(n,\lceil r-n^{\gamma}\rceil)\ge 1/(2e)^2-{\lambda \over n^{2\gamma-1}}.
\end{align*}
Since $n\ge n_1+n_2\ge 2n_0$, applying the bound (\ref{eq:lambda-order-1}) on $\lambda$, this is at least $1/(3e)^2$ for $n_0$ large enough. Then applying
the concentration bound (\ref{eq:Chebyshev-tail}), for $n_0$ large enough
$\E{H_n} \ge r-n^{\gamma}-2n^{1-\delta}\ge r-3n^{1-\delta}$. Applying Lemma~\ref{lemma:EHn-close-lambda} we conclude
$\lambda_n n \ge  r-5n^{1-\delta}$ provided $n_0$ is large enough. Dividing by $n$ we obtain for sufficiently large $n_0$
\begin{align*}
\lambda_n \ge \lambda - 5n^{-\delta}=\min(\lambda_{n_1},\lambda_{n_2})-5 n^{-\delta},
\end{align*}
concluding the proof.  
\end{proof}

Next we lower bound $\lambda_{n+1}$ in terms of $\lambda_n$.
\begin{lemma}\label{lemma:from-n-to-n+1}
There exists $n_0$ such that for all  $n\ge n_0$
\begin{align*}
\lambda_{n+1} \ge \lambda_n-4(n+1)^{-\delta}.
\end{align*}
\end{lemma}

\begin{proof}
Let $D_{n+1}$ be the number of clauses containing variable $x_{n+1}$ in $\Phi(n+1,m)$. For $m\le 2^{k+1}n$ we have 
$\inf_{n}\pr(D_{n+1}=0)\triangleq \alpha>0$. In fact 
$D_{n+1}$ has a Poisson distribution in the limit with $n$ independent parameter, but this will not be explicitly needed.  
We have $\pr(\Phi(n+1,m)\in {\rm SAT} | D_{n+1}=0)=\pr(\Phi(n,m)\in {\rm SAT})$.
Thus for every $m\le 2^{k+1}n$
\begin{align*}
\pr(H_{n+1}\ge m) & = \pr(\Phi(n+1,m)\in {\rm SAT})\\
& \ge \pr(\Phi(n+1,m)\in {\rm SAT} | D_{n+1}=0)\alpha \\
&= \pr(\Phi(n,m)\in {\rm SAT})\alpha \\
& \ge \pr(H_n\ge m)\alpha.
\end{align*}
Applying this to $m=\lambda_n n$
\begin{align*}
\pr(H_{n+1}\ge \lambda_n n) &\ge \pr(H_n\ge \lambda_n n)\alpha \\
&=\alpha/2.
\end{align*}
We claim that  for $n_0$ sufficiently large and all $n\ge n_0$, $\E{H_{n+1}}\ge \lambda_n n-(n+1)^{1-\delta}$. Indeed, otherwise
\begin{align*}
\alpha/2\le 
\pr(H_{n+1}\ge \lambda_n n)
\le 
\pr(H_{n+1}\ge \E{H_{n+1}}+(n+1)^{1-\delta})\le 1/(n+1)^\eta,
\end{align*}
which is violated when $n\ge n_0$ and $n_0$ is large enough. Combining with Lemma~\ref{lemma:EHn-close-lambda}
we obtain
\begin{align*}
\lambda_{n+1}(n+1)\ge \lambda_n n-3(n+1)^{1-\delta}.
\end{align*}
Dividing by $n+1$ and using $\lambda_n/(n+1)\le (n+1)^{-\delta}$, which holds for sufficiently large $n_0$ and by (\ref{eq:lambda-order-1}), 
 we obtain the result.
\end{proof}

\subsection{From superadditivity to the existence of a limit of normalized medians $\lambda_n$. Proof of the main theorem}
We complete the proof of the main theorem in this section by establishing the existence of a limit $\lim_n\lambda_n$ of normalized medians using 
the near superadditivity property.  The following variant  of the superadditivity property is adopted from Abbe and Montanari~\cite{abbe2014concentration}

\begin{prop}\label{prop:superadditive}
Suppose $a_n, n\ge 1$ is a non-negative bounded sequence of real values satisfying the following properties. 
There exists $n_0$ and $\delta>0$ such that for all $n_1,n_2\ge n_0$
\begin{align}
a_{n_1+n_2}\ge \min(a_{n_1},a_{n_2})-5 (n_1+n_2)^{-\delta}, \label{eq:a-superadd}
\end{align}
and for all $n\ge n_0$
\begin{align}
a_{n+1} \ge a_n-4(n+1)^{-\delta}. \label{eq:an+1 vs an}
\end{align}
Then for all $r\ge n_0$ and all $n\ge r^2$
\begin{align}
a_n\ge a_r-{14r^{-\delta}\over 1-(2/3)^\delta}. \label{eq:n-vs-r}
\end{align}
As a result  the limit $a^*\triangleq \lim_{n\to\infty} a_n$ exists. 
\end{prop}

We first show how this claim implies the main theorem.
\begin{proof}[Proof of Theorem~\ref{theorem:main-result}]
Corollary~\ref{coro:superadditive} and Lemma~\ref{lemma:from-n-to-n+1} imply that our sequence $\lambda_n$
satisfies the assumptions of the Proposition~\ref{prop:superadditive} with $\delta=(k-2)/(4k)$ chosen as in Proposition~\ref{prop:variance}.
Thus the limit $\lambda^*=\lim_n \lambda_n$ exists.  By (\ref{eq:lambda-order-1}) we have $0<\lambda^*\le 2^{k}$.
Fix any $\epsilon>0$ such that $\lambda^*(1+\epsilon)<2^{k+1}$.  
By Lemma~\ref{lemma:EHn-close-lambda} we can find large enough $n_0$ such that for all $n\ge n_0$
$\E{H_n}\le \lambda_n n+2n^{1-\delta}$. Again taking $n_0$ sufficiently large we obtain for all $n\ge n_0$
\begin{align*}
\E{H_n}\le  \lambda^*(1+\epsilon/2)n.
\end{align*}
Take $m=\lceil \lambda^*(1+\epsilon)n \rceil\le 2^{k+1}n$. 
Then applying the concentration bound Proposition~\ref{prop:variance} for large enough $n_0$ and all $n\ge n_0$
\begin{align*}
\pr(\Phi(n,m)\in {\rm SAT}) &= \pr(H_n\ge m) \\
&\le
\pr(H_n-\E{H_n}\ge \lambda^*(\epsilon/2)n)\le Cn^{-\eta}.
\end{align*}
Here we used the uniform bound (\ref{eq:lambda-order-1}) on $\lambda_n$.
Thus $\lim_{n\to\infty}\pr(\Phi(n,m)\in {\rm SAT})=0$ and the one side of the assertion of the theorem is established.

The other side is established similarly. Taking $n_0$ sufficiently large we obtain for all $n\ge n_0$
\begin{align*}
\E{H_n}\ge  \lambda^*(1-\epsilon/2)n.
\end{align*}
We set $m=\lceil \lambda^*(1-\epsilon)n \rceil$. Then
\begin{align*}
\pr(\Phi(n,m)\notin {\rm SAT}) &= \pr(H_n<m) \\
&\le
\pr(H_n-\E{H_n}\le -\lambda^*(\epsilon/2)n)\le Cn^{-\eta}.
\end{align*}
Thus $\lim_{n\to\infty}\pr(\Phi(n,m)\in {\rm SAT})=1$, and the other side of the assertion of the theorem is established.
\end{proof}

Thus it remains to  prove Proposition~\ref{prop:superadditive}. 

\begin{proof}[Proof of Proposition~\ref{prop:superadditive}]
The existence of the limit $\lim_n a_n$ is immediate from (\ref{eq:n-vs-r}). Thus we focus on proving (\ref{eq:n-vs-r}).

Given any positive integer $R$  consider a binary tree $\mathcal{T}(R)$ constructed recursively as follows. 
Every node $u$ of this tree is marked with a positive integer denoted by  $R_u$. The root $o$ of the tree $\mathcal{T}(R)$ is marked
with value $R$ itself.
When $R=1$ the associated tree is a single node $u$,
with $R_u=R=1$. 
Assuming the tree is defined for all $1\le R'<R$, we construct it for $R$ recursively as follows. 
The root of the  tree $\mathcal{T}(R)$ has two children $o_1, o_2$. 
Suppose $R=2\hat R$ is even. Then the  subtrees emanating from both $o_j, j=1,2$ are $\mathcal{T}(R')$. 
In particular, $R_{o_j}=R'$. On the other hand, when $R=2R'+1$ is odd, the two children $o_j, j=1,2$ of the root $o$ of $\mathcal{T}(R)$ 
are roots of subtrees $\mathcal{T}(R')$ and $\mathcal{T}(R'+1)$, respectively. In particular, the two children $o_1$ and $o_2$ are marked
with $R'$ and $R'+1$, respectively. As a result, the leaves of every tree $\mathcal{T}(R)$ are marked with value $1$.
Since in either even or odd case we have $(1/3)R\le R'\le (2/3)R$, the depth $T$ of the tree is logarithmic in $R$,
though this fact  will not be required for the remainder of the proof.

Now, fix any $r\ge n_0$ and $n\ge r^2$. Finding  $q\ge r, 0\le t\le r-1$ such that $n=qr+t$,  
we can write $n=(q-t)r+t(r+1)$. Note $q-t>0$.  
We assume  $t>0$. The case $t=0$ is treated similarly, and the proof is omitted.
Applying (\ref{eq:a-superadd}) with $n_1=(q-t)r$ and $n_2=t(r+1)$ we have 
\begin{align}
a_n\ge \min(a_{(q-t)r},a_{t(r+1)})-5n^{-\delta}. \label{eq:a-n-vs-q-t-r}
\end{align}
Now consider the two  binary trees $\mathcal{T}(q-t)$ and $\mathcal{T}(t)$ defined as above. 
Starting from the left tree $\mathcal{T}(q-t)$ with the root   marked with $q-t$, apply (\ref{eq:a-superadd})  recursively down the  tree until 
the leaves of the left tree are reached. 
For example, in the first step, assuming $q-t$ is even we obtain
\begin{align*}
a_{(q-t)r}\ge a_{{q-t\over 2}r}-5((q-t)r)^{-\delta}.
\end{align*}
Assuming, on the other hand,  $q-t$ is odd, in the first step  we obtain  instead
\begin{align}
a_{(q-t)r}\ge \min\left(a_{{q-t-1\over 2}r},a_{{q-t+1\over 2}r}\right)-5((q-t)r)^{-\delta}. \label{eq:top-error}
\end{align}
Iterating we obtain
\begin{align*}
a_{(q-t)r}\ge a_r-\max_P E_P
\end{align*}
where $E_P$ is the additive correction  accumulated across a  path $P$ from the root to a leaf. 
 From the relation $R'/R\le 2/3$ between a mark of a node and the
mark of any of its child, that we observed above, every next correction term up the path 
is at most $(2/3)^\delta$ times the previous correction term. The largest correction is associated with leaves and is $5r^{-\delta}$. 
Thus we obtain
\begin{align*}
E_P\le 5r^{-\delta}\sum_{\ell \ge 0}(2/3)^{\delta \ell} = {5r^{-\delta}\over 1-(2/3)^\delta}.
\end{align*}
We obtain 
\begin{align*}
a_{(q-t)r}\ge a_r- {5r^{-\delta}\over 1-(2/3)^\delta}.
\end{align*}
Similarly, we obtain 
\begin{align*}
a_{t(r+1)} \ge a_{r+1}-{5(r+1)^{-\delta}\over 1-(2/3)^\delta}.
\end{align*}
Applying (\ref{eq:an+1 vs an}) we obtain 
\begin{align*}
a_{t(r+1)} \ge a_r-{5(r+1)^{-\delta}\over 1-(2/3)^\delta}-4(r+1)^{-\delta}.
\end{align*}
Applying this to (\ref{eq:a-n-vs-q-t-r}) we obtain
\begin{align}
a_n &\ge a_r-{5r^{-\delta}\over 1-(2/3)^\delta}-4(r+1)^{-\delta}-5n^{-\delta} \\
&\ge
a_r-{14r^{-\delta}\over 1-(2/3)^\delta},
\end{align}
completing the proof. 
\end{proof}

\section{Some open questions}\label{section:Open-problems}
Several questions appear to be within the reach of the current techniques, in particular ones  developed 
in~\cite{GCarenini} and~\cite{k-sat-OpenAI}. In particular, it seems very plausible that the approach extends to other models of random
constraint satisfaction problems such as proper coloring of a graph. 

Furthermore, while the validity of the limiting satisfiability threshold
$\lambda_k^*$ is now established, it is not known whether it coincides with ground state energy limit, the existence of which was established
in~\cite{BayatiGamarnikTetali}. Specifically, it is shown in~\cite{BayatiGamarnikTetali}  that for every fixed value $d$, the minimum number of violated clauses
in the random formula $\Phi(n,\lceil dn\rceil)$ converges to some limit $\gamma_k(d)$ when normalized by $n$. Define
$d^*_k$ as $\sup\{d: \gamma_k(d)=0\}$.  $d^*$ is strictly positive since  $\gamma_k(d)$ is zero with high probability  when $d<\lambda_k^*$. 
What remains unclear, though is extremely plausible, is whether $d^*=\lambda_k^*$. Specifically, at the present stage, 
we cannot rule out a range of values $d$ 
above $\lambda_k^*$ for which on the one hand the formula is not satisfiable with high probability, 
but on the other hand solutions exist which violate only sublinearly in $n$ many clauses also with high probability.
Provided an  exponential concentration of the hitting time $H_n$ around its expectation holds, such a conclusion would be an immediate
corollary of Theorem~\ref{theorem:main-result} and the log-superadditivity property Theorem~\ref{theorem:BGT-log-additivity} from~\cite{BayatiGamarnikTetali}.
Unfortunately, the current concentration bound from Proposition~\ref{prop:variance} is too weak for this purpose. 

Finally, unfortunately not all random constraint satisfaction problems such as random k-SAT and spin glasses are amenable to the interpolation technique.
The key challenge is that the convexity property such as one in (\ref{eq:convexity-1}) and (\ref{eq:convexity-2}) does not always hold. A notable counter-example
is the Minimum Bisection problem in  sparse random graph, see~\cite{Huang2018} for the definition and some partial results. 
Specifically, given a sparse \ER graph $\G(n,d/n)$ let ${\rm MB}(n,d)$ be its minimum bisection value. The expectation of this value is $O(n)$ trivially
and it is natural to conjecture that the limit $\lim_n \E{{\rm MB}(n,d)}/n$ exists. 
Unfortunately,  the counterparts of the expressions (\ref{eq:convexity-1}) and (\ref{eq:convexity-2}) for the Minimum Bisection problem do not exhibit the
convexity property. As a result the validity of the conjecture remains open.

\bibliographystyle{amsplain}
\bibliography{bibliography-10.2026}

\end{document}